\documentclass[a4paper]{article}
\usepackage{geometry}                
\usepackage[parfill]{parskip}    

 \usepackage[usenames]{color}
\usepackage{amsthm}
\usepackage{amssymb}
\usepackage{latexsym}
\usepackage{euscript}
\usepackage{amsmath}
\usepackage{amstext}
\usepackage{amsgen}
\usepackage{amsbsy}
\usepackage{amsopn}
\usepackage{amsfonts}
\usepackage[sans]{dsfont}
\usepackage[all]{xy}
\usepackage{hyperref}
\usepackage{graphicx}
\usepackage{subfigure}
\usepackage{epstopdf}
\title{Transversal Homotopy Monoids of Stratified Complex Projective Space}
\author{Conor Smyth}

\newtheorem{thm}{Theorem}

\newtheorem{lem}[thm]{Lemma}
\newtheorem{prop}[thm]{Proposition}

\newtheorem{remark}[thm]{Remark}
\newtheorem{cor}[thm]{Corollary}

\newtheorem{defn}[thm]{Definition}

\begin{document}
\maketitle
\section{Introduction}
This paper, which follows from the work in \cite{Wlf}, will give a geometric description of $\psi_n(\Bbb{CP}^k)$, the $n$th transversal homotopy monoid of $k$-dimensional complex projective space, where $\Bbb{CP}^k$ is stratified as follows, $\Bbb{CP}^0\subset\Bbb{CP}^1\subset\ldots\subset\Bbb{CP}^k.$ 
Transversal homotopy monoids are defined as classes of based transversal maps into Whitney stratified spaces up to equivalence through such maps. 
In a homotopy group each class has an inverse obtained by precomposing a representative with a reflection of the sphere. 
This is not the case in transversal homotopy. To see this consider a class in the fundamental group represented by a based path which transversally crosses a codimension $1$ stratum. We can choose the dual of this class represented by the same path but traversed in the opposite direction and compose the two. This is not homotopic through transversal maps to the identity because we would be required to pull the path back across the stratum and it would at some point fail to be transverse. And so transversal homotopy in some sense detects the strata within a space. The invariants we obtain are not in general groups but monoids with duals.

These transversal homotopy monoids do not easily admit either an algebraic or geometric description. For instance if we take $\Bbb{S}^n$ to be the $n$-sphere stratified by a point and its complement, and look at $\psi_n(\Bbb{S}^n)$, the $n$th transversal homotopy monoid, we can describe it algebraically as the free (commutative for $n\geq 2$) monoid with duals on one generator. They can also be viewed geometrically as isotopy classes of framed codimension $n$ submanifolds on the $n$-sphere. See \cite{Wlf} for details. 
On the other hand $\psi_{n+1}(\Bbb{S}^n)$ can be described as framed codimension $n$ submanifolds inside $n+1$ space. For instance $\psi_2(\Bbb{S}^1)$, that is classes of based transversal maps of $S^2$ into $\Bbb{S}^1$ (up to homotopy through such maps), can be represented as nested circles in $\Bbb{R}^2$ up to ambient isotopy; but an algebraic description is not so apparent. Another interesting example is $\psi_3(\Bbb{S}^2)$, the set of all framed links up to isotopy.
\footnotetext{This paper was completed under the close supervision of Jon Woolf who I thank for his generosity of time. The work is funded by the Leverhulme 
Trust (grant reference: F/00 025/AH).}
$\Bbb{CP}^k$ admits a natural stratification as described above, with depth $\geq 1$ and contractible strata. By depth we mean the longest chain of strata contained within a space. It is also an interesting way of generalising $\Bbb S^2$ i.e.  $\Bbb{CP}^1\cong \Bbb S^2$, and so we would like to find a nice geometric description of its transversal homotopy monoids $\psi_n(\Bbb{CP}^k)$ by relating transversal maps of the sphere into $\Bbb {CP}^k$ with the pullback stratifications of the sphere.
In doing so we must answer the following question:
What are the necessary and sufficient conditions for a stratification of $S^n$ to be the pullback stratification of a transversal map $f:S^n\to\Bbb{CP}^k$?

\section{Preliminaries} 
\subsection{Some Necessary Background}

A \emph{Whitney stratification of a manifold} $X$ is a disjoint decomposition $X$ = $\bigcup_{\alpha}X_{\alpha}$ of $X$ into submanifolds (which are not necessarily connected, but which must each have a fixed dimension for all of their connected components) that satisfies the following four axioms:

\begin{enumerate}
\item \emph{Local finiteness}. The decomposition is locally finite, i.e. every point $x \in X$ has a neighbourhood \emph{U} with the property that $U\cap X_{\alpha}$ is empty for all but a finite number of strata $X_{\alpha}$.

\item \emph{The axiom of the frontier}. If one stratum $X_\alpha$ has a non-empty intersection with the closure $\overline{X_{\beta}}$ of another stratum $X_\beta$, then $X_\alpha$ lies entirely within  $\overline{X_\beta}$.

\item \emph{Whitney's condition A}. Suppose $X_{\alpha}$ lies entirely in the closure of $X_{\beta}$. \\ Suppose that $x_1, x_2, x_3, \dots$ is a sequence of points in $X_\beta$ which converges to a point $y$ in $X_\alpha$. Suppose further that the sequence of tangent spaces $T_{x_1}X_\beta, T_{x_2} X_\beta, T_{x_3} X_\beta \dots$ converges as a sequence of subspaces of the tangent space TX to X, to a ``limiting tangent space"  $\tau \subseteq T_y X$. Then the tangent space $T_{y} X_ \alpha$ is contained in the limiting space $\tau$.

\item \emph{Whitney's condition B}. Suppose that $X_\alpha$ lies in the closure of $X_\beta$. Suppose $x_1,x_2,x_3\dots$ is a sequence of points in $X_{\beta}$ which converge to a point $y$ in $X_\alpha$, and $y_1,y_2,y_3,\dots$ is a sequence of points in $X_\beta$ which also converge to $y$. Suppose as before the sequence of tangent spaces $T_{x_1}X_\beta, T_{x_2} X_\beta, T_{x_3} X_\beta ,\dots$ converges, as a sequence of subspaces of the tangent space TX to X, to  $\tau \subseteq T_y X$. Suppose further that the sequence of secant lines $\overline{x_1y_1},\overline{ x_2y_2},\overline{ x_3y_3},\dots$ converges to a limiting line $l\subseteq T_yX$. Then the limiting line is contained in the limiting tangent space $\tau$.
\end{enumerate}
In fact Whitney's condition $B$ implies condition $A$ but condition $A$ is used more in applications.

Two smooth submanifolds $A$ and $B$ of a smooth manifold $M$ are said to intersect \emph{transversely} if for any point $x\in A\cap B$ we have, $$T_x A + T_x B = T_x M.$$ The sum need not be direct, i.e. we do not require $$T_x A \cap T_x B = \{0\}.$$
If the intersection is empty we say that $A$ and $B$ are vacuously transverse.
More generally, smooth maps $f:A\to M$ and $g:B\to M$ are transverse if for each $p\in A$ and $q\in B$ where $f(p)=g(q)$ we have
$$df(T_p A)+dg(T_q B)=T_{f(p)=g(q)}M.$$
This is equivalent to the composite $$T_p A\stackrel{df}{\longrightarrow}T_{f(p)}M\to T_{f(p)}M/dg(T_q B)$$ being surjective.
Transversality of submanifolds corresponds to the case where $f$ is the inclusion of $A$ and $g$ is the inclusion of $B$.

A map $f:X\to W$ where $X$ is a manifold and $W$ a Whitney stratified manifold, is transversal if and only if for each stratum $S\subset W$ the restriction $$f|_s:S\to X$$ is transverse to each stratum of $X$.

\begin{thm}\label{trans}If the smooth map $f:M\to M'$ of manifolds is transversal to a submanifold $A\subset M'$, then the preimage $f^{-1}(A)$ is a submanifold of $M$. Moreover the codimension of $f^{-1}(A)$ in $M$ equals the codimension of $A$ in $M'$. Furthermore we can identify $Nf^{-1}A\cong f^*NA$ where $NA$ and $Nf^{-1}A$ denote the normal bundles to $A$ in $M'$ and to $f^{-1}A$ in $M$ respectively.
\end{thm}
\begin{defn}
A \emph{tubular neighborhood} of a submanifold $A$ in $M$ is an embedding of the normal bundle  $NA$ of $A$ into $M$, i.e., $f: NA\to M$, where the image of the zero section of the normal bundle is equal to $A$ in $M$. 
\end{defn}
\begin {thm}[Tubular Neighbourhood Theorem]\label{tube}
Every embedded  submanifold of $\Bbb R^n$ has a tubular neighbourhood. 
\end{thm}
\begin{remark}\label{rem1}
Given any manifold $M$ with a submanifold $A$, we can choose a metric on $M$. Then $NA\cong(TA)^{\perp}\subset M$, is a subbundle of $TM$. The exponential map of the metric provides a tubular neighbourhood. More precisely we can choose $\epsilon >0$ such that the open $\epsilon$-disk bundle in $NA$ is smoothly embedded in $M$ by the exponential map. Therefore a choice of metric and suitable $\epsilon$ gives a well-defined tubular neighbourhood, $NA\to U_A.$
Since the space of such choices is path-connected, it follows that any two such tubular neighbourhoods $$M\xleftarrow{\alpha}NA\xrightarrow{\beta}M$$ are connected by a family $\alpha_t:NA\to M$ of tubular neighbourhoods such that        
$\alpha_0=\alpha$ and $\alpha_1=\beta.$
\end{remark}

\begin{thm}[Whitney's Approximation Theorem for Manifolds]\label{Approx}
Let $M$ and $M'$ be smooth manifolds and let $F:M\to M'$ be a continuous map. Then $F$ is homotopic to a smooth map $G:M\to M'$. If $F$ is smooth on a closed subset $A\subset M$ the homotopy can be taken relative to $A$. 
\end{thm}

Proof of Theorem \ref{trans} can be found for example in \cite{Guil}, page 28, and of Theorems $3$ and $5$ in \cite{Lee}, pages 255 and 259, respectively.

\subsection{Transversal Homotopy Monoids}
Let $X$ be a Whitney stratified manifold and fix a generic basepoint $x\in X$, i.e. a base point in an open stratum. For $n\in \Bbb N$ we define the \emph{$n^{th}$ transversal homotopy monoid} $\psi_{n}(X)$ to be the set of transversal maps $[0,1]^n \to X$ which map some neighbourhood of the boundary to the basepoint $x$ up to the equivalence relation generated by homotopy through transversal maps. We denote the class of a transversal map $f$ by $[f]$.

For $n\geq 1$ juxtaposition gives $\psi_n(X)$ the structure of a monoid:
$[f]\cdot[g]=[f\cdot g]$ where

 \[f\cdot g : [0,1]^n\to X:
 \left \{ \begin{array}{ll} 
 f(t_1,\dots,2t_n) & t_n\in[0,\frac{1}{2}]\\
g(t_1,\dots,2t_n -1) & \textrm{$t_n\in(\frac{1}{2},1]$.} \end{array}\right.\]

Insisting that a neighbourhood of the boundary maps to $x$ ensures that the juxtaposition is smooth; it is clearly transversal. The class of the constant map to $x$ is the unit for juxtaposition. For $n\geq 2$ we could choose to juxtapose in any of the coordinate directions. However the usual Eckmann-Hilton argument shows that these monoidal structures are all the same and are commutative.

Remark: Since any continuous map $f:X\to Y$, where $X$ and $Y$ are manifolds, is homotopic to a smooth one, (see for example \cite{Lee}, page 257, Theorem 10.21), if $X$ has only one stratum then $$\psi_n(X)\cong \pi_n(X).$$

\subsection{The Pullback Stratification}
A transversal map $f : [0, 1]^n\to X$ determines an induced stratiÞcation of $[0, 1]^n$ by the subsets $f^{-1}S$ where $S\subset X$ is a stratum.
Each stratum $A$ of this induced stratification maps to a stratum in $X$, which we denote by $S_A$. The 
restriction of $f$ to $A$ determines a homotopy class $[f |_A] \in [A, S_A]$ of maps from 
$A$ to $S_A$. Furthermore, using the fact that $f$ is transverse to $S$, 
the derivative of $f $ induces an isomorphism of bundles 
$$(f |_A)^{*} N S_A \cong NA $$
where $ NS_A$ is the normal bundle of the stratum $S_A$ in $X$. Up to 
isomorphism the pullback bundle on the left depends only on $[f 
|_A] \in[A, S_A]$. With respect to the induced stratification the map $f$ is stratified, that is, it sends strata to strata.

\begin{prop} \label{prop} Suppose $h : W 
\times [0, 1] \to X$ is a homotopy through transversal
maps. Then there is a (continuous) ambient isotopy of $W 
\times [0, 1]$ starting at 
the identity and ending at a stratum-preserving homeomorphism from the stratification induced by $h$ to the product of the stratification on $W$ induced from 
$h_0 : W 
\to X$ and the trivial stratification on $[0, 1]$. 
\end{prop}
This is a restatement of Thom's First Isotopy Lemma, proved by Mather in the notes \cite{TFIL}.
\\
 \begin{lem} (see \cite{Wlf} section $3.3.$)
 An element $[f]\in \psi_n(X)$ where $n\geq 1$ is invertible if and only if the stratification induced by $f$ is trivial.
 \begin{proof} It follows from Proposition ~\ref{prop} that if the induced stratification of a transversal map $f:[0,1]^n\to X$ is not trivial then the induced stratification of any other representative is non-trivial too. Hence the condition is invariant under homotopies through transversal maps. Furthermore if the stratification induced by $f$ is non-trivial then so is that induced by any composite $f\cdot g$ and so $[f]$ cannot be invertible.
 Conversely if the stratification induced by $f$ is trivial then $f$ maps $[0,1]^n$ into the open stratum containing the basepoint and the usual inverse of homotopy theory provides an inverse in $\psi_n(X)$.
\end{proof}
\end{lem} 
So $\psi_n(X)$ is not in general a group rather it is a dagger monoid, that is, a unital associative monoid with involution $a\mapsto a^{\dag}$ such that $1^{\dag}=1, (ab)^{\dag}=b^{\dag}a^{\dag}$.
The involution on $\psi_n(X)$ for $n>0$ is given by $[f]^\dag=[f^\dag]$ where \\ $$f^\dag:[0,1]^n\to X:(t_1,\cdots ,t_n)\mapsto f(t_1,\cdots,t_{n-1},1-t_n).$$

\section{Construction} \label{Con}
A class in $\psi_n(\Bbb{CP}^k)$ is represented by a transversal map $f: S^n \to\Bbb{CP}^k.$
The pre-image stratification arises from a filtration of the form $X_0 \subset \cdots \subset X_k=S^n$ where
$$X_i=f^{-1}\Bbb{CP}^i$$
and the strata are the differences $X_i - X_{i+1}$. 

What are the necessary conditions on such a filtration to arise this way?
Firstly, $f^{-1}(\Bbb{CP}^{i-1})$ must be a manifold of real codimension $2$ inside $f^{-1}(\Bbb{CP}^{i}).$
This is due to Theorem \ref{trans} in which we take the submanifold to be $\Bbb{CP}^i$.

Next, since the normal bundle to $\Bbb{CP}^i$ inside $\Bbb{CP}^{i+1}$ is a complex line bundle, we can classify it (up to isomorphism) by its first Chern class.
For an oriented submanifold $A^i\subset M^n$ we write $[A]$ for the Poincar\`e dual class to $A$ in $H^{n-1}(M;\Bbb Z).$

Using naturality to pull back the Chern class along $f$, and then the fact that $f$ is a transversal map we have,
\begin{eqnarray*}
f^*c_1(N_{\Bbb{CP}^{i+1}}\Bbb{CP}^i)&=&c_1(f^*N_{\Bbb{CP}^{i+1}}\Bbb{CP}^i)\\
&=&c_1(N_{f^{-1}\Bbb{CP}^{i+1}}f^{-1}\Bbb{CP}^i)\\
&=&c_1(N_{X_{i+1}}X_i)
\end{eqnarray*}

But it is well-known that $c_1(N_{\Bbb{CP}^{i+1}}\Bbb{CP}^i)=[\Bbb{CP}^{i-1}]\in H^2(\Bbb{CP}^i;\Bbb Z)$ so,
\begin{eqnarray*}
f^*c_1(N_{\Bbb{CP}^{i+1}}\Bbb{CP}^i)&=&f^*[\Bbb{CP}^{i-1}]\\
&=&[f^{-1}\Bbb{CP}^{i-1}]\\
&=&[X_{i-1}]\in H^2(f^{-1}\Bbb{CP}^i;\Bbb Z)=H^2(X_i;\Bbb Z).
\end{eqnarray*}
Again we have used the fact that $f$ is transversal to identify $f^*[\Bbb{CP}^{i-1}]$ with $[f^{-1}\Bbb{CP}^{i-1}]$.
So we see that for each $X_i$ it is necessary that $c_1(N_{X_{i+1}}X_i)=[X_{i-1}]\in H^2(X_i;\Bbb Z)$ and that $X_i$ must be of real codimension $2$ inside $X_{i+1}$.

For the remainder of this paper we will use the fact that the first Chern class $c_1$, of a complex line bundle over a manifold can be identified with the Euler class $e$ of its underlying oriented real bundle. 
\begin{prop}\label{Prop}
Suppose $f|_{X_i}:X_i\to\Bbb{CP}^i$ is defined and transversal, with the property that the pullback stratification is obtained from the filtration $$X_0 \subset \cdots \subset X_i,$$ and that $X_j$, $0\leq j\leq i+1$, is a codimension $2n-2j$ submanifold and $N_{X_{i+1}}X_i)$ can be oriented so that $e(N_{X_{i+1}}X_i)=[X_{i-1}]\in H^2(X_i;\Bbb Z).$
Then we can extend to $f:X_{i+1}\to \Bbb{CP}^{i+1}$. Moreover given a choice of homotopy class of an isomorphism, $N_{X_i+1}X_i\cong f^* N_{\Bbb{CP}^{i+1}} \Bbb{CP}^i$, the extension is unique up to homotopy through transversal maps. 

\begin{proof}
To construct such a map we use the idea of a collapse map from the Pontryagin--Thom construction. 
We define the collapse map on $X_{i+1}$ as follows. 
Let $U_{X_{i}}\subset X_{i+1}$ be a tubular neighbourhood of $X_i$ inside $X_{i+1}$. We must choose an isomorphism, $$X_{i}\subset U_{X_{i}}\cong N_{X_{i+1}}X_{i}.$$ Next since 
\begin{eqnarray*}
e(N_{X_{i+1}}X_i)&=&[X_{i-1}]\\
&=&[f^{-1}\Bbb{CP}^{i-1}]\\
&=& f^*[\Bbb{CP}^{i-1}]\in H^2(\Bbb{CP}^i;\Bbb Z)\\
&=&f^*e(N_{\Bbb{CP}^{i+1}}\Bbb{CP}^{i})\\
&=&e(f^*N_{\Bbb{CP}^{i+1}}\Bbb{CP}^i)\\
\end{eqnarray*}
we can choose an isomorphism,
$$N_{X_{i+1}}X_{i}\cong f^* N_{\Bbb{CP}^{i+1}}\Bbb{CP}^{i}.$$ We can now extend $f$ to a map $U_{X_{i}}\cong N_{\Bbb{CP}_{i+1}}\Bbb{C}_{i}$ by defining 
$$(x,w)\mapsto(f(x),w)$$ 
where $w$ is in the fibre of $f^*N_{\Bbb{CP}^{i+1}}\Bbb{CP}^i$ at $x$ which by the definition of a pullback bundle is canonically the same as the fibre of $N_{\Bbb{CP}^{i+1}}\Bbb{CP}^i$ at $f(x)$. From this we can complete the following, $$U_{X_i}\xrightarrow{\textrm{tub. neighbourhood}}N_{X_{i+1}}X_i\xrightarrow{\textrm{choice}} f^*N\Bbb{CP}^i\xrightarrow{\textrm{canonical}} N_{\Bbb{CP}^{i+1}} \Bbb{CP}^i\xrightarrow{\textrm{tub. neighbourhood}}\Bbb{CP}^{i+1}.$$

Let  $D_i\subset U_i$ be the image of a closed disk bundle inside $N_{X_{i+1}}X_i$.
Consider $f|_{\partial D_i}:\partial D_i\to \Bbb{CP}^{i+1}\setminus\Bbb{CP}^i\cong \Bbb C^{i+1}.$ Since $\Bbb C^{i+1}$ is contractible we can extend this to a continuous map $$f:X_{i+1}\setminus \textrm{int}D_i\to  \Bbb{CP}^{i+1}\setminus\Bbb{CP}^i.$$
Moreover, such extensions are unique up to homotopy.
This map is smooth in a neighbourhood of $X_i$ and continuous outside so by the Whitney Approximation Theorem it is homotopic relative to $U_{X_i}$ to a smooth map.

It was necessary to make some choices in the construction. We needed to choose a tubular neighbourhood of $X_i$ inside $X_{i+1}$ but by Remark \ref{rem1} we see that any two choices of tubular neighbourhood are connected by a family of tubular neighbourhoods so our choice is unique up to homotopy. 
We made a similar choice of tubular neighbourhood of $\Bbb{CP}^{i}$ in $\Bbb{CP}^{i+1}$ but again this choice was unique up to homotopy using the same reasoning.
We also needed to choose an isomorphism $N_{X_i+1}X_i\cong f^* N_{\Bbb{CP}^{i+1}} \Bbb{CP}^i$ but so long as we fix the homotopy class, the resulting extention is unique up to homotopy through transversal maps.
\end{proof}
\end{prop}
We now show that the choice of homotopy class of the isomorphism $N_{X_i+1}X_i\cong f^* N_{\Bbb{CP}^{i+1}} \Bbb{CP}^i$ is equivalent to choosing a framing of $X_i\times\Bbb{R}^2$. This is analogous to the choice of framing required in the Pontryagin--Thom construction.

\begin{lem}\label{Lem} The set of isomorphisms (up to homotopy) $$N_{X_i+1}X_i\cong f^* N_{\Bbb{CP}^{i+1}} \Bbb{CP}^i$$ is (non-canonically) the set of framings of $X_i\times\Bbb R^2$, i.e. bundle isomorphisms of $X_i\times\Bbb R^2.$
\begin{proof} Choose one isomorphism $\phi:N_{X_i+1}X_i\to f^* N_{\Bbb{CP}^{i+1}} \Bbb{CP}^i$ and let $L$ be a complex line bundle on $X_i$, where $c_1(L)=-c_1(f^*N_{\Bbb{CP}^{i+1}}\Bbb{CP}^i$). There are bijections,
\begin{displaymath}
\xymatrix{\theta\ar@{|->}[d]&\in &\textrm{Isom}(N_{X_i+1}X_i, f^* N_{\Bbb{CP}^{i+1}} \Bbb{CP}^i)\ar@{<=>}^\cong[d]\\
\theta\circ\phi^{-1}&\in&\textrm{Isom}(f^* N_{\Bbb{CP}^{i+1}} \Bbb{CP}^i,N_{X_i+1}X_i)\ar@{<=>}[d]\\
&&\textrm{Isom}(f^* N_{\Bbb{CP}^{i+1}} \Bbb{CP}^i\otimes L,f^* N_{\Bbb{CP}^{i+1}} \Bbb{CP}^i\otimes L)\ar@{<=>}[d]\\
&&\textrm{Isom}(X_i\times\Bbb C,X_i\times\Bbb C)
}
\end{displaymath}
This identification respects homotopies.
\end{proof}
\end{lem}
\begin{cor} Given initial data and choices (as above) there is a transversal map $$S^n\xrightarrow{f}\Bbb{CP}^k$$ such that the pullback stratification is obtained from the filtration $$X_0 \subset \cdots \subset X_k=S^n,$$ where $X_i$ is a codimension $2k-2i$ submanifold and $e(N_{X_{i+1}}X_i)=[X_{i-1}]\in H^2(X_i;\Bbb Z).$ 
\begin{proof}
On $X_0$, $f$ is unique (up to homotopy) and maps all of $X_0$ to a point. 
Note that $N_{X_1} X_0\cong f^*N_{\Bbb{CP}^1} \Bbb{CP}^0$ must be trivial. The map is therefore well-defined on a neighbourhood of $X_0$.
From here we proceed inductively using the method described in Proposition \ref{Prop}.
\end{proof}
\end{cor}

\section{Conclusion}
We have shown that each class in $\psi_n(\Bbb{CP}^k)$ is represented by such a filtration, and every such filtration represents some class.
It remains to consider the relation on such filtrations corresponding to homotopy through transversal maps. This relation is smooth isotopy, see \cite{Wlf} (Proposition 2.13 and Remark 2.14). Conversely if we have an isotopy $$a:S^n\times [0,1]\to S^n$$ and a filtration $X_0\subset\cdots\subset X_k$ of $S^n$ (of required kind) then there is a filtration $$\alpha^{-1}X_0\subset\cdots\subset\alpha^{-1}X_k$$ of $S^n\times[0,1].$
We can apply the construction of Section \ref{Con} to this filtration of $S^n\times[0,1]$ to obtain a homotopy through transversal maps $S^n\times[0,1]\to \Bbb{CP}^k.$
So $\psi_n(\Bbb{CP}^k)$
has a geometric description as smooth isotopy classes of filtrations of $S^n$ of the form $X_0 \subset \cdots \subset X_k,$ where $X_i$ is a codimension $2k-2i$ submanifold such that the Euler class $e(N_{X_{i+1}}X_i)=[X_{i-1}]\in H^2(X_i;\Bbb Z).$

\begin{thebibliography}{Woo10}

\bibitem[GP74]{Guil}
Victor Guillemin and Alan Pollack.
\newblock {\em Differential topology}.
\newblock Prentice-Hall Inc., Englewood Cliffs, N.J., 1974.

\bibitem[Lee03]{Lee}
John~M. Lee.
\newblock {\em Introduction to smooth manifolds}, volume 218 of {\em Graduate
  Texts in Mathematics}.
\newblock Springer-Verlag, New York, 2003.

\bibitem[Mat70]{TFIL}
J.~Mather.
\newblock Notes on topological stability,
  1970.
\newblock {\em Available from the following website {\tt www.math.princeton.edu/facultypapers/mather/}}

\bibitem[Woo10]{Wlf}
Jonathan Woolf.
\newblock Transversal homotopy theory.
\newblock {\em Theory Appl. Categ.}, 24:No. 7, 148--178, 2010.

\end{thebibliography}

\end{document}